\documentclass[11pt,a4paper]{article}
\usepackage{amssymb,epsfig}

\def\cz{{\mathbb C}}

\def\gz{{\mathbb Z}}
\def\pr{{\mathbb P}}

\def\skp{\hspace{1pt}}
\def\eps{\varepsilon}
\def\ph{\varphi}

\def\id{\mathop{\rm Id\skp}}

\def\MC{{\mathop{\rm MC}}}
\def\HMC{{\mathop{\rm HMC}}}
\def\||{|\!|}
\newcommand{\diagl}[1]%
  {\makebox[0cm]{${\scriptstyle#1\ }
  \hspace{-5pt}\downarrow\hspace{-5pt}
  \phantom{\scriptstyle#1\ }$}}
\newcommand{\diagr}[1]%
  {\makebox[0cm]{$\phantom{\ \scriptstyle#1}
  \hspace{-5pt}\downarrow\hspace{-5pt}{\ \scriptstyle#1}$}}
\newcommand{\updiagl}[1]%
  {\makebox[0cm]{${\scriptstyle#1\ }
  \hspace{-5pt}\uparrow\hspace{-5pt}
  \phantom{\scriptstyle#1\ }$}}
\newcommand{\updiagr}[1]%
  {\makebox[0cm]{$\phantom{\ \scriptstyle#1}
  \hspace{-5pt}\uparrow\hspace{-5pt}
  {\ \scriptstyle#1}$}}

\newcommand{\cala}{{\cal A}}

\newcommand{\calf}{{\cal F}}

\newcommand{\calo}{{\cal O}}

\newtheorem{proposition}{Proposition}[section]
\newtheorem{theorem}[proposition]{Theorem}
\newtheorem{lemma}[proposition]{Lemma}
\newtheorem{corollary}[proposition]{Corollary}
\newtheorem{remark}[proposition]{Remark}

\newtheorem{definition}[proposition]{Definition}
\newcommand{\pf}{{\em Proof. }}
\newcommand{\qed}{{{\hfill$\diamond$}\vspace{1.5ex}}}

\title{On hyperelliptic $C^\infty$-Lefschetz fibrations of four-manifolds}
\author{Bernd Siebert and Gang Tian}
\date{February 28, 1999}
\begin{document}
\maketitle
\begin{abstract}
    We show that hyperelliptic symplectic Lefschetz fibrations
    are symplectically birational to two-fold covers of rational
    ruled surfaces, branched in a symplectically embedded surface.
    This reduces the classification of genus 2 fibrations
    to the classification of certain symplectic submanifolds in
    rational ruled surfaces.
\end{abstract}

\tableofcontents
\vspace{4ex}


\addcontentsline{toc}{section}{Introduction}
\noindent
{\Large\bf Introduction}
\vspace{1.5ex}

\noindent
A \emph{differentiable Lefschetz fibration (DLF)} of an oriented,
differentiable four manifold $M$ is a proper differentiable map
$q:M\to S^2$ with finitely many critical points $Q_1,\ldots, Q_\mu\in M$
in disjoint fibers and such that locally near $Q_r$ there exist complex
coordinates $z,w$ on $M$ and $t$ on $S^2$ with $q:(z,w)\mapsto zw$.
If $(z,w)$ and $t$ can be chosen compatible with the orientations
of $M$ and $S^2$, then $q:M\to S^2$ is called a \emph{symplectic
Lefschetz fibration (SLF)}. According to an observation of Gompf
a DLF is a SLF iff there exists a symplectic form
$\omega$ on $M$ compatible with $q$ in the sense that smooth
fibers of $q$ are symplectic submanifolds with respect to
$\omega$. An isomorphism of two DLF's or SLF's $q:M\to S^2$ and
$q':M'\to S^2$ is an oriented diffeomorphism of total spaces
$\Psi:M\to M'$ respecting the fibration structures.

Given two SLF's (or DLF's) $q:M\to S^2$, $q':M'\to S^2$ and an
(orientation preserving) diffeomorphism of smooth fibers
$\Phi: q^{-1}(s_0)\simeq {q'}^{-1}(s'_0)$ the {\em fiber connected
sum} of $q$ and $q'$ is the SLF (DLF) obtained by identifying the
boundaries of $q^{-1}(S^2\setminus B_\eps(s_0))$ and
${q'}^{-1}(S^2\setminus B_\eps(s'_0))$ via $\Phi$.
The origin of this paper is a question of one of the authors (G.T.)
if every SLF can be decomposed into a fiber sum of algebraic
Lefschetz fibrations. This is true for elliptic fibrations ($g=1$,
\cite{moishezon}). To the authors knowledge there
still exists no example of a SLF that is known to not decompose
in this way. Discussions during the 1993/1994 stay of the
first named author at the Courant Institute suggested
that at least for fibrations by curves of genus 2, one should
expect algebraicity when all singular fibers are irreducible.
The work took on more concrete form during the academic year
1997/1998, where we looked at this problem from the point of
view of symplectic geometry inside $S^2$-bundles $P\to S^2$.
This paper contains the first, rather elementary step of
reduction to the classification problem of symplectic submanifolds of
$P$ up to isomorphism. The authors have a program to attack
the slightly stronger problem of classification up to isotopy
with the technique of pseudo-holomorphic curves.
\vspace{2ex}

\noindent
The results of this paper were presented at the G\"okova
Geometry-Topology conference (May 1998) and at the
Geometry conference at Schloss Ringberg (June 1998).

While this paper was under preparation we learnt that
Ivan Smith proved similar results in his Oxford thesis.
We apologize for any duplication of results.


\section{Hyperelliptic fibrations and branched covers}
Our main focus in this paper will be on so-called hyperelliptic SLF's,
which are defined in terms of the monodromy as we will now describe.
Let $s_1\ldots,s_\mu\in S^2$ be the critical values of $q$. The
restriction of $q$ to $S^2\setminus \{s_1,\ldots,s_\mu\}$ is a fiber
bundle  with fibers a closed oriented surface $\Sigma$ of some genus
$g$,  the \emph{genus} of the DLF or SLF. The corresponding monodromy
representation
\[
  \rho:\ \pi_1(S^2\setminus\{s_1,\ldots,s_\mu\})\
  \longrightarrow\ \MC_g=\mathrm{Diff}^+{(\Sigma)}/
  \mathrm{Diff}^0(\Sigma)
\]
sends closed loops running counterclockwise once around some $s_r$
({\em simple loops}) to a Dehn-Twist. Since for loops $\gamma,\gamma'$
the product $\gamma\gamma'$ means running through $\gamma$ first,
we let the mapping class group act on $\Sigma$ from the right to make
$\rho$ a homomorphism. Note also that we do admit Dehn-twists
along contractible loops. These occur iff an irreducible component of the
singular fiber is an embedded sphere of self-intersection $-1$.
The Dehn-twist is right-handed if the fibration is oriented
locally near $s_r$, otherwise left-handed. The monodromy
representation determines the fibration in the following precise way.
\begin{theorem}\label{thm.kas}
    \textnormal{\cite{kas}}\ \ There is a one-to-one
    correspondence between the set of isomorphism classes of
    SLF's of genus $g\ge2$ and with $\mu$ singular fibers and the set of
    representations $\rho: \pi_1 (S^2\setminus\{s_1,\ldots,s_\mu\})
    \to \MC_g$ sending simple loops to right-handed Dehn-twists
    modulo composition with inner automorphisms of $\MC_g$.
    The same result holds for $g=1$ provided there is at least
    one irreducible singular fiber.
\end{theorem}
By the same proof the analogous statement holds for DLF's and
arbitrary Dehn-twists. Note also that $\pi_1
(S^2\setminus\{s_1,\ldots,s_\mu\})$ can be generated by $\mu$
simple loops $\gamma_1,\ldots,\gamma_\mu$ obeying the single
relation $\gamma_1\cdot\ldots \cdot\gamma_\mu=1$, so a $DLF$
(SLF) can be defined simply by a relation between (right-handed)
Dehn-twists.

To define the hyperelliptic mapping class group we fix once and for all a
two-sheeted branched cover $\kappa:\Sigma\to  S^2$, necessarily
having $2g+2$ branch points. Then the \emph{hyperelliptic
mapping class group} $\HMC_g\subset \MC_g$ is the subgroup
generated by diffeomorphisms of $\Sigma$ that are compatible
with $\kappa$, that is, which commute with the hyperelliptic involution
swapping the sheets of $\kappa$. Finally a DLF or SLF is
\emph{hyperelliptic} if the monodromy representation is equivalent
(by an inner automorphism of $\MC_g$) to one taking values in $\HMC_g$.

The hyperelliptic mapping class group is a $\gz/2\gz$-extension of the
mapping class group $\MC(S^2,2g+2)$ of $S^2$ with $2g+2$ marked points.
The copy of $\gz/2\gz$ is of course generated by the hyperelliptic
involution, while the map to $\MC(S^2,2g+2)$ is the map to the
induced diffeomorphism of the base, mapping the set of $2g+2$
branch points to itself (well-definedness of this map requires some
work \cite{birman2}). To write down a presentation of
$\MC(S^2,2g+2)$ let $x_1,\ldots,
x_{2g+1}$ be a standard generating set of Dehn-twists of
$(S^2=\cz\cup\{\infty\};P_1,\ldots,P_{2g+2})$, $P_r=\exp(\frac{2\pi i
r}{2g+2})$, represented by a diffeomorphism $\Phi_r$
that is identical outside of a small neighbourhood of the straight line
connecting $P_r$ and $P_{r+1}$, while these points get exchanged. Then
$\MC(S^2;2g+2)$ has the presentation \cite[p.164]{birman}
\[
  \MC(S^2,2g+2)\ =\ \Bigg\langle x_1,\ldots, x_{2g+1}
  \Bigg| \begin{array}{c}x_i x_j=x_j x_i\ \ {\scriptstyle |i-j|\ge2}\,;\
           x_i x_{i+1}x_i=x_{i+1}x_i x_{i+1}\\
           x_1\ldots x_{2g+1}x_{2g+1}\ldots x_1=1\\
           (x_1\ldots x_{2g+1})^{2g+2}=1\end{array}\Bigg\rangle
\]
The $x_i$ lift to right-handed Dehn-twists on $\Sigma$, that we will denote
by the same notation. These Dehn-twists still fulfill the braid
relations, but $x_1\ldots x_{2g+1}x_{2g+1}\ldots x_1$ becomes the
hyperelliptic involution $I$. In fact, $\HMC_g$ has the presentation
\cite{birman2}
\[
  \HMC_g\ =\ \Bigg\langle x_1,\ldots, x_{2g+1}, I
  \Bigg| \begin{array}{c}x_i x_j=x_j x_i\ \ {\scriptstyle |i-j|\ge2}\,;\
           x_i x_{i+1}x_i=x_{i+1}x_i x_{i+1}\\
           I=x_1\ldots x_{2g+1}x_{2g+1}\ldots x_1;\
	   I x_i=x_i I;\ I^2=1\\
           (x_1\ldots x_{2g+1})^{2g+2}=1\end{array}\Bigg\rangle
\]
This should be contrasted to the complexity of $\MC_g$ ($g\ge3$), which needs
one more Dehn-twist for a generating set and a number of more involved
relations
\cite{wajnryb}. Our main suggestion in this
paper is to study hyperelliptic SLF's by means of
branch loci in $S^2$-fibrations over $S^2$. This is based on
Theorem~\ref{red.to.branch} below.

Before stating the theorem we need a few more notions. A map $\sigma:
\tilde M\to M$ between oriented, differentiable four-manifolds will be
called the
{\em blow-up of $M$ in $P\in M$} if $\sigma$ is a diffeomorphism away from
$P$ and if there is a neighbourhood $U$ of $P$  such that
$\sigma^{-1}(U)\to U$
is oriented diffeomorphic to the blow up of a point in an open set in
$\cz^2$. The
contracted set $\sigma^{-1}(P)$ will be called {\em exceptional divisor}.

It is a well-known fact that holomorphic fibrations over $\pr^1$ with
general fiber $\pr^1$ and smooth total space are simply
blow ups of $\pr^1$-bundles over $\pr^1$. As an ad-hoc
definition we define $S^2$-fibrations over $S^2$
in the differentiable setting as blow ups of $S^2$-bundles over
$S^2$. One class of $S^2$-fibrations will be of special interest, namely
those where each singular fiber is {\em of type $\Lambda_2$}, which by
definition means that locally it is obtained from an $S^2$-bundle by blowing
up one point and then another general point on the exceptional divisor thus
obtained. Such fibers consist of a chain of embedded spheres with
self-intersections  $-1$, $-2$, $-1$.
\begin{theorem}\label{red.to.branch}
    A SLF $q:M\to S^2$ is hyperelliptic iff it fits into a
    commutative diagram of the form
    \begin{eqnarray}\begin{array}{rccccc}
      \tilde q:&\tilde M&\stackrel{\tilde \kappa}{\longrightarrow}&\tilde P&
      \stackrel{\tilde p}{\longrightarrow}&S^2\\
      &\diagl{\sigma} &&\diagr{\rho}&&\diagr{\id}\\
      q:&M&\stackrel{\kappa}{\longrightarrow}&P&
      \stackrel{p}{\longrightarrow}&S^2
    \end{array}\label{diagthm}\end{eqnarray}
    with
    \begin{enumerate}
    \item
      $\sigma$ is the blow up of $M$ in the nodes of {\em reducible} singular
      fibers of $q$
    \item
      $\tilde p:\tilde P\to S^2$ is an $S^2$-fibration with all
      singular fibers of type $\Lambda_2$
    \item
      $\tilde \kappa$ is a two-fold branched cover with branch locus
      $B+B_E\subset\tilde P$ where ${\tilde \kappa}^{-1}(B_E)$ is the
      exceptional divisor of $\sigma$ (so $B_E$ is a disjoint union
      of $(-2)$-spheres in the fibers of $\tilde p$) and $B$ has no
      components in the fibers of $\tilde p$
    \item
      $\rho$ contracts $B_E$, leading to an orbifold $P$ with an
      $A_1$-singularity (locally isomorphic to $\cz^2$ modulo the
      diagonal $\gz/2\gz$-action) for each component of $B_E$
    \end{enumerate}
\end{theorem}
Before we can give the proof we need some preparations. We begin with
the braid group $B(S^2,2g+2)$ of $S^2$ with $2g+2$ strands
\cite[\S1.5]{birman}. It has the presentation
\[
   B(S^2,2g+2)\ =\ \Bigg\langle x_1,\ldots, x_{2g+1}
  \Bigg| \begin{array}{c}x_i x_j=x_j x_i\ \ \ {\scriptstyle |i-j|\ge2}\,;\
           x_i x_{i+1}x_i=x_{i+1}x_i x_{i+1}\\
           x_1\ldots x_{2g+1}x_{2g+1}\ldots x_1=1
  \end{array}\Bigg\rangle\,.
\]
Our notations are consistent in that the element $x_i$
represents the movement that exchanges the $i$-th and
$(i+1)$-th strands by a positive half-twist. So the natural
homomorphism $B(S^2,2g+2)\to\MC(S^2,2g+2)$ induces the
presentation of $\MC(S^2,2g+2)$ above.
One can show that $(x_1\cdot\ldots\cdot x_{2g+1})^{2(2g+2)}=1$ in this
group \cite[p154f]{birman}. Geometrically this relation can be seen by
pulling $2g+2$ parallel strands across $\infty$: This operation
results in two full twists of the bunch of strands (in fact,
$(x_1\cdot\ldots\cdot
x_{2g+1})^{2g+2}$ represents one twist). On the other hand, on the
side of the mapping class group $\MC(S^2,2g+2)$ a full twist induces
a rotation of $S^2$ by $2\pi$, and hence is trivial.
So $B(S^2,2g+2)$ is another (central) extension of $\MC(S^2,2g+2)$
by $\gz/2\gz$, non-isomorphic to $\HMC_g$.

The importance of $B(S^2,2g+2)$ in our situation comes from the fact
that homomorphisms
\[
  \rho':\ \pi_1(\cz\setminus\{s_1,\ldots,s_\mu\})\ \longrightarrow\
  B(S^2,n)
\]
(modulo inner automorphisms of $B(S^2,n)$) classify isotopy classes
of closed submanifolds $B^*\subset S^2\times (\cz\setminus
\{s_1,\ldots,s_\mu\})$ mapping to $\cz\setminus
\{s_1,\ldots,s_\mu\}$ as $n$-sheeted unbranched covers via monodromy.
We will use the following lemma to relate the monodromy of the branch
locus to the monodromy of the fibration.
\begin{lemma}\label{compatibility.of.mon}
    For $\eps<\min\{\mathrm{dist}_{r\neq r'}(s_r,s_{r'})\}/2$ put
    $S^*=\cz\setminus B_\eps(\{s_1,\ldots,s_\mu\})$, $P^*=S^*\times
    S^2$ and let $\kappa: M^*\to P^*$ be a two-to-one cover with branch
    locus $B^*\subset P^*$. We furthermore assume that the composition
    $B^*\hookrightarrow P^*\to S^*$ is a $(2g+2)$-sheeted unbranched cover.
Let
    $p:P^*\to S^*$ and $q:M^*\to S^*$ be the corresponding fiber
    bundles. Choose $s_0\in S^*$ and compatible isomorphisms
    $q^{-1}(s_0)\simeq \Sigma$ (the closed surface of genus $g$ from
    above) and $(p^{-1}(s_0),p^{-1}(s_0)\cap B^*)\simeq
    (S^2;P_1,\ldots,P_{2g+2})$.

    Then the monodromy representations describing $M^*\to S^*$ and
    $B^*\hookrightarrow P^*\to S^*$
    \begin{eqnarray*}
        \rho:\ \pi_1(S^*) & \longrightarrow & \HMC_g  \\
        \rho':\ \pi_1(S^*) & \longrightarrow & B(S^2,2g+2)
    \end{eqnarray*}
   are compatible: They induce the same homomorphism into
   $MC(S^2,2g+2)$.
\end{lemma}
\pf
Let $\gamma: S^1\to S^*$ be a closed loop, $\gamma(0)=s_0$. Let
$\calf$ be a horizontal foliation for $p_\gamma:
P_\gamma:=S^1\times_\gamma P^*\to S^1$ such that $B_\gamma=
S^1\times_\gamma B^*\subset P_\gamma$ is a leaf. The monodromy
$\rho'(\gamma)$ is represented by the flow of $S^2$ obtained by
following the leaves and using the given trivialization
$P_\gamma=S^1\times S^2$. The corresponding element of
$\MC(S^2,2g+2)$ can then be represented by the self-diffeomorphism
$\Phi_{\gamma,\calf}$ from lifting the identical path $S^1\to S^1$.
Note that $\Phi_{\gamma,\calf}$ maps $B\cap p^{-1}(s_0)$ onto itself
because $B_\gamma$ is a leaf.

On the other hand, the fact that $B_\gamma$ is a leaf allows us to
lift $\calf$ to a foliation $\tilde\calf$ of $M^*$. Thus
$\rho(\gamma)$ is represented by a self-diffeomorphism
$\Psi_{\gamma,\tilde\calf}$ of $\Sigma=q^{-1}(s_0)$ commuting with
$\Phi_{\gamma,\calf}$ via $\kappa$:
\[\begin{array}{ccc}
    \Sigma & \stackrel{\Psi_{\gamma,\tilde \calf}}{\longrightarrow} &
\Sigma  \\
    \diagl{\kappa} &  & \diagr{\kappa}  \\
    S^2 & \stackrel{\Phi_{\gamma,\calf}}{\longrightarrow} & S^2
\end{array}\]
This is what we claimed.
\qed
\vspace{2ex}

\noindent
Another ingredient in the proof of the theorem is a list of
local models for (possibly singular) branch loci
$B\subset \Delta\times S^2$ producing all
occurring Dehn-twists. For Lefschetz fibrations the monodromy
can be represented by a Dehn-twist $\Psi$ along an embedded circle
$C\simeq S^1\subset\Sigma$, the vanishing cycle, cf.\ e.g.\
\cite[Thm.2.1]{kas}. Up to a diffeomorphism of $\Sigma$ there are
$2+[\frac{g}{2}]$ different kinds of embedded circles. They can be
distinguished
by checking if $\Sigma\setminus C$ is connected or disconnected
and in the latter case by the genera of the two components.
Accordingly there are only $2+[\frac{g}{2}]$ essentially different
kinds of Dehn-twists. For monodromies of Lefschetz pencils along
simple loops the type of Dehn -twist can be read off from the genera of
the irreducible components of the enclosed singular fiber.
If $\Psi$ commutes with the hyperelliptic
involution $I$ (the hyperelliptic case) then $C$ is isotopic to $I(C)$
and we can even assume $C=I(C)$:
\begin{lemma}\label{comp.van.cycles}
    Let $\kappa:\Sigma\to S^2$ be a two-fold branched cover and
    $I:\Sigma\to \Sigma$ the hyperelliptic involution exchanging the
    branches of $\kappa$. Assume that $\gamma: C\subset \Sigma$ is a
    simple closed curve which is homotopic to $I(C)$. Then $C$ is isotopic
    to a simple closed curve $C'\subset \Sigma$ with $C'=I(C')$.
\end{lemma}
\pf
We endow $\Sigma$ with the complex structure making
$\kappa$ a holomorphic covering, so $\Sigma$ becomes
a hyperellipic Riemann surface and $\kappa$ the 1-canonical map.
The Poincar\'e metric $h$ on $\Sigma$ is then invariant under
the hyperelliptic involution $I$. Recall the classical fact that
in any free homotopy class of closed loops on a closed surface
with negative curvature there is a unique closed geodesic loop.
We define $C'\subset \Sigma$ as the unique geodesic loop with
respect to $h$ freely homotopic to $C$. By $I$-invariance of $h$
the loop $I(C')$ is also a geodesic loop, homotopic to $I(C)$.
But by hypothesis, $I(C)$ is homotopic to $C$, and hence $C'=I(C')$.
The proof is finished with the theorem of Baer saying that simple
closed loops are isotopic if and only if they are homotopic \cite{baer}.
\qed
\vspace{1ex}

Let us therefore assume $C=I(C)$. We distinguish two cases: (1)
$C\cap  \mathrm{Fix}(I)=
\emptyset$\ \ (2) $C\cap\mathrm{Fix}(I)\neq\emptyset$. We will see
that the former case covers exactly the disconnecting circles.

In the first instance, $\bar C=\kappa(C)$ is an embedded circle in
$S^2\setminus\{P_1,\ldots, P_{2g+2}\}$ and $\kappa: C\to \bar C$ is a
two-fold unbranched cover. The simple Dehn-twist $\Psi$ along $C$
thus pushes down to a {\em double} Dehn-twist $\Phi$ along $\bar C$:
$\kappa\circ \Psi=\Phi\circ\kappa$. Moreover, since the monodromy
($\in\gz/2\gz$) of $\kappa$ along $\bar C$ is non-trivial, $\bar C$
divides $\{P_1,\ldots,P_{2g+2}\}$ into two subset of {\em odd}
cardinality, say $\{P_1,\ldots,P_{2h+1}\}$,
$\{P_{2h+2},\ldots,P_{2g+2}\}$, $h\ge 0$.

In the second instance, write $\kappa: z\mapsto w=z^2$ in appropriate
complex coordinates $z=x+iy$, $w$ centered at some $P\in
C\cap\mathrm{Fix} (I)$ and in $\kappa(P)$ respectively, so
$I:z\mapsto -z$. Possibly after a rotation we may write $C$ locally
as graph of a real function $\ph$ with $\ph(0)=0$, that is
$C=\{x+i\ph(x)\}$. The property $C=I(C)$ then implies
$\ph(-x)=-\ph(x)$ and, locally,
\[
  \bar C\ =\ \{(x^2-\ph^2(x))+2ix\ph(x)\}\,.
\]
This shows that $\bar C$ is an embedded, compact, connected curve in
$S^2$ with one end for each point $P_i\in \bar C$. Since $\bar C$ has
thus at least one end, the classification of plane, embedded curves
implies that $\bar C$ is a closed intervall connecting exactly two
points $P_i\neq P_j$. It is a well-known, elementary fact that a Dehn-twist
$\Phi$ along a circle enclosing $\bar C$ (say the boundary of
an $\eps$-neighbourhood of $\bar C$) commutes with a diffeomorphism
$\Psi'$ of $\Sigma$ that is isotopic to $\Psi$.
We have thus described all Dehn-twists in $\HMC_g$ as lifts of
explicit diffeomorphisms $\Phi$ of $(S^2; P_1,\ldots,P_{2g+2})$.

We now write down elements of the braid group $B(S^2,2g+2)$ producing
these diffeomorphisms $\Phi$ up to isotopy. To be
precise we write $S^2=\cz\cup\{\infty\}$ and assume
\[
  \bar C=\{|z|=2\}\,;\quad P_r=e^\frac{2\pi r i}{2h+1}
  \mathrm{\ \ for\ }r=1,\ldots,2h+1\,;\quad
  |P_r|>1\mathrm{\ for\ }r>2h+1
\]
in the first instance and
\[
  \bar C = [-1,1]\,;\quad P_1=-1\,;\quad P_2=1\,;\quad
  |P_r|>1\mathrm{\ for\ }r>2
\]
($i=1,j=2$)
in the second instance. The double Dehn-twist along $\bar C$ (case 1)
is then induced by the braid ($t\in[0,1]$)
\[
  P_r(t)=e^{(\frac{r}{2h+1}+2t)\cdot2\pi i}
  \mathrm{\ \ for\ }r=1,\ldots,2h+1\,;\quad
  P_r(t)=P_r\mathrm{\ for\ }r>2h+1\,,
\]
while the Dehn-twist along $\bar C$ (case 2) is induced by the half-twist
\[
  P_1(t)=-e^{\pi i t}\,;\quad P_2(t)=e^{\pi i t}\,;\quad
  P_r(t)=P_r\mathrm{\ for\ }r>2\,.
\]
These braids occur as monodromies of the following holomorphic
curves $\bar B\subset (2\Delta)\times\pr^1$ along $\gamma(t)=e^{2\pi i
t}\in 2\Delta$:
\begin{eqnarray}\label{standard.curves}
    \begin{array}{ll}
        (B_{h,g-h})\hspace{10ex}&\begin{array}{c} \Big\{(s,z)\in
                        (2\Delta)\times \pr^1\,\Big|\,
	\prod_{r=1}^{2h+1}(z- e^\frac{2\pi r i}{2h+1}s^2)=0\Big\}\\[1.5ex]
	\cup(2\Delta)\times\{P_{2h+2,}\ldots,P_{2g+2}\}\end{array}\\[4ex]
        (B_\mathrm{irr})&\begin{array}{c}\Big\{(s,z)\in (2\Delta)\times
                        \pr^1\,\Big|\, z^2-s=0\Big\}\\[1.5ex]
	\cup (2\Delta)\times\{P_1,\ldots,
	\hat P_i,\ldots,\hat P_j,\ldots,P_{2g+2}\}\end{array}
    \end{array}
\end{eqnarray}
where in the second line the entries with a hat are to be
omitted. The second curve is smooth, but unless $h=0$ the first curve
has a singularity of multiplicity $2h+1$
that in algebraic geometry is traditionally called {\em infinitely
close multiple point} (here: of multiplicity $2h+1$). After one blow
up the strict transform of this curve becomes an ordinary singular
point of the same multiplicity $2h+1$ (locally the union of $2h+1$ smooth
curves
with {\em disjoint} tangent planes), so another blow up desingularizes
(hence {\em infinitely} close). The first non-trivial case $g=2$,
$h=1$ is depicted in Figure~1 below.\\
\epsfxsize=9cm\hspace*{2.15cm}\epsfbox{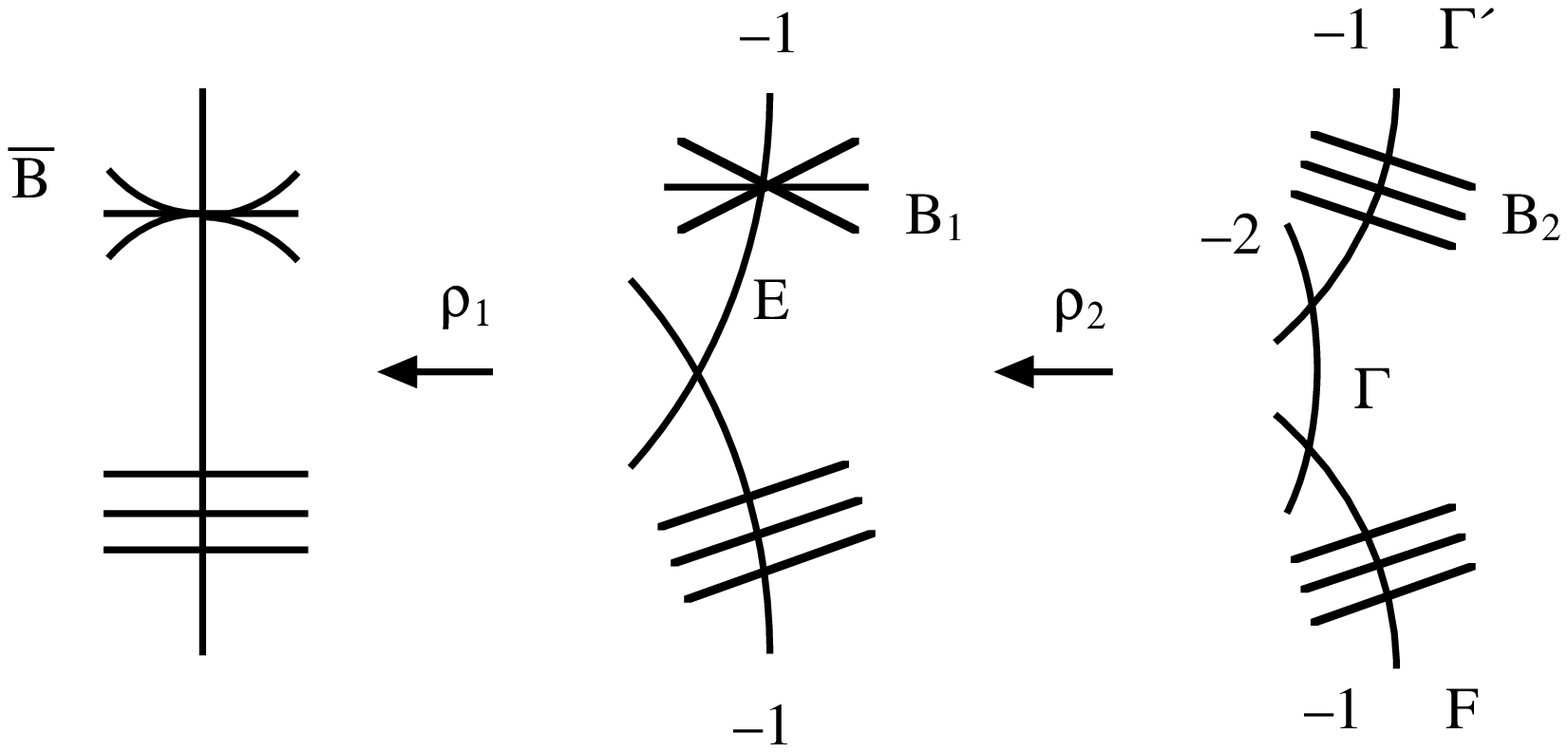}
\\[-3cm]
\mbox{}\hfill Fig.~1\\[3ex]
The numbers are the self-intersection numbers of the various
components. This should be compared to statement (3) of the
proposition.
\vspace{2ex}

As third preparatory remark let $p:P=\pr(\calo\oplus \calo(n))\to
\pr^1=S^2$ be a $\pr^1$-bundle and $B\subset P$ a pseudomanifold
projecting to $S^2$ as unbranched cover (of degree $d$) away from
finitely many points. We want to relate the homology class of $B$
to a ``virtual number of critical values'' of $B\to S^2$ that we are
going to define momentarily. In the sequel we will reserve $H,S_\infty,
F\subset P$ for (the classes of) a positive section, a negative section and a
fiber of $P$ respectively, so $H^2=n$ and $S_\infty=-n$.
Since $H_2(P,\gz)=\langle H,F\rangle$ we may write $B\sim d\cdot H+l\cdot F$
in homology. To determine $l$
we describe $B$ by its {\em developing map}
\[
  \Phi_B:\ \pr^1\ \longrightarrow\ S^d_{\pr^1} P\,,
\]
sending $s\in S^2=\pr^1$ to the 0-cycle $p^{-1}(s)\cap B$ (with
multiplicities).
Here $S^d_{\pr^1} P = P\times_{\pr^1}\ldots\times_{\pr^1}P/S_d$
is the $d$-fold fibered symmetric product of $P$
with itself. Since $S^d\pr^1\simeq\pr^d$ the target of $\Phi_B$ is a
$\pr^d$-bundle over $\pr^1$. Its homology is thus generated by a fiber class
$\eta$ and some horizontal class $\xi$, for which we take the reduced divisor
of points in $S^d_{\pr^1} P$ parametrizing (zero-dimensional) cycles with
support
meeting $S_\infty\subset P$. Let $D\subset
S^d_{\pr^1} P$ be the discriminant locus corresponding to cycles
with higher multiplicities. Note that $D$ is the
push-forward of the generalized diagonal in $P\times_{\pr^1}\ldots
\times_{\pr^1}P$ under the symmetrization map, so $D$ is
an irreducible and reduced divisor. We define the {\em virtual number
of critical values} of $B\to S^2$ by
\[
  \mu_\mathrm{virt} (B)\ :=\ \deg \Phi_B^*D\,.
\]
\begin{lemma}\label{lem.birt.branch.pts}
    $D\sim 2(d-1)\xi+nd(d-1)\eta$.
\end{lemma}
\pf
We start writing $D\sim a\xi+b\eta$ for some $a,b\in \gz$.
To determine $a$ and $b$ we consider two special cases.
The first one is $B\sim dH$ represented as union of graphs of $d$ pairwise
different
multiples $\lambda_i s$ of a section $s$ (i.e.\ viewing $P$ as
compactification of
$\calo_{\pr^1}(n)$\,). The discriminant of this $d$-tuple is a non-zero
constant (the
discriminant of the $\lambda_i$) times $s^{d(d-1)}$, which
counted with multiplicities has $nd(d-1)$ zeroes.
This shows $b=nd(d-1)$. The other case is $B\sim (d-1)H+S_\infty\sim dH-nF$
represented by only $d-1$ of the sections and $S_\infty$. We obtain
$-an+b=n(d-1)(d-2)$, so $a=2(d-1)$, which can also be seen as the degree
of the discriminant $D\subset S^d\pr^1\simeq\pr^d$.
\qed

Since $\Phi_B$ is a section of $S^d_{\pr^1}P\to \pr^1$ we have $\deg
\Phi_B^*\eta =1$. On the other hand, $\deg\Phi_B^* \xi=S_\infty\cdot
B= l$ by the definition of $\xi$. We conclude
\begin{proposition}\label{numbersing}
  \nopagebreak
  The number $\mu_\mathrm{virt}$ of virtual branch points of the
  projection $B\subset\pr(\calo\oplus\calo(n))\rightarrow\pr^1$
  is related to the homology class $B=dH+lF$ in the following way\\[1.3ex]
  \mbox{}\hfill$\mu_\mathrm{virt}\ =\ 2(d-1)l+nd(d-1)\,.$\hfill$\diamond$
\end{proposition}
\vspace{1ex}

For later reference let us also record quickly the contribution to
$\mu_\mathrm{virt}(B)$ of our local standard models: For the case
$B_\mathrm{irr}$ of smooth critical points one easily checks in
local coordinates of $S^d\pr^1$ that the contribution is $+1$. The
case $B_\mathit{h,g-h}$ with a point of multiplicity $2h+1$
is not so obvious. We may however reduce to the previous case by a
holomorphic deformation of $B$ to a smooth curve. To this end we first
translate each of the $2h+1$ smooth irreducible components slightly to
make them meet transversally. Since each pair of irreducible
components meets in 2 points we obtain
\[
  2\cdot{{2h+1} \choose 2}\ =\ 2h(2h+1)
\]
nodes (ordinary double points). In turn each node can be smoothed,
leading to 2 simple critical points of the projection to the base (that is,
of type
$B_\mathrm{irr}$). The contribution to $\mu_\mathrm{virt}(B)$ is
therefore
\[
  4h(2h+1)\,.
\]
Note that this number coincides with the number of generators in the
corresponding braid monodromy $(x_1\cdot\ldots\cdot x_{2h})^{2(2h+1)}$.
In fact, this deformation argument is another way to deduce this
monodromy (cf.\ \cite[p.83ff]{chakiris} for $h=1$).
\vspace{2ex}

\noindent
We are now ready to prove the theorem.
\vspace{1ex}

\noindent
{\em Proof of Theorem~\ref{red.to.branch}.}\ \
The case $g=1$, which would require some extra considerations, follows
from the classification of Moishezon and Livn\'e \cite[II.2]{moishezon},
so we assume $g\ge2$. We arrange the critical values $s_1,\ldots,s_\mu$ of
$q:M\to  S^2$ to be the $\mu$-th roots of unity ($S^2=\cz\cup\{\infty\}$). For
sufficiently small $\eps>0$ let $S^*=\cz\setminus
B_\eps(\{s_1,\ldots,s_\mu\})$ and
$\rho:\pi_1(S^2\setminus\{s_1,\ldots, s_\mu\})\to \HMC_g$ the
monodromy representation. Let $\gamma_r$ be the loop that runs from
$0$ once counterclockwise around $s_r$ and which does not enclose any
$s_{r'}$, $r'\neq r$. Then $\pi_1(S^*)$ is the free group generated by
$\gamma_1,\ldots,\gamma_\mu$, while $\pi_1(S^2\setminus\{s_1,\ldots,
s_\mu\})$ is the quotient by the relation $\gamma_1\cdot\ldots\cdot
\gamma_\mu=1$. The monodromy
representation is then given by Dehn-twists $\tau_r=\rho(\gamma_r)\in
\HMC_g$ fulfilling
\[
  \tau_1\cdot\ldots\cdot\tau_\mu\ =\ 1\,.
\]
Let $\bar\tau_r$ be the image of $\tau_r$ in $\MC(S^2,2g+2)$. There are
distinguished lifts $\hat\tau_r\in B(S^2,2g+2)$ as follows: By the discussion
above we know that a non-trivial $\tau_r$ is conjugate either to $x_1$ or to
$(x_1\cdot\ldots\cdot x_{2h})^{4h+2}$ ($1\le h\le [g/2]$), and the same
expression holds for $\bar\tau_r$. Now simply represent the conjugating
element as word in the $x_r$, and interpret the whole word as element of
$B(S^2,2g+2)$. While the lift of the conjugating element is not unique,
the product is well-defined because $B(S^2,2g+2)$ is
a {\em central} extension of $\MC(S^2,2g+2)$.

By this process we might loose triviality of the monodromy along the
circle $\gamma_1\cdot\ldots\cdot\gamma_\mu$ at infinity, but
at least
\[
  \hat\tau_1\cdot\ldots\cdot \hat\tau_\mu\ \in\ \{1,(x_1\ldots
  x_{2g+1})^{2g+2}\}\,.
\]
The $\hat\tau_r$ determine a $(2g+2)$-fold unbranched cover
$B^*\hookrightarrow S^*\times S^2\to S^*$ having $\hat\tau_r$ as
braid monodromy. In fact, let $D=\{(z_1,\ldots,z_{g+2})\in
(S^2)^{2g+2}\,|\, z_r=z_j$ for some $i\neq j\}$ be the generalized
diagonal in $(S^2)^{2g+2}$, and
\[
  A_r=(A_r^k)_k:\ [0,1]\ \longrightarrow\ (S^2)^{2g+2}\setminus D
\]
with  $A_r(0)=(P_1,\ldots,P_{2g+2})$, $A_r(1)=(P_{\sigma_r(1)},\ldots,
P_{\sigma_r(2g+2)})$ be a movement of the $2g+2$
branch points of $\kappa: \Sigma\to S^2$
representing $\hat\tau_r$. So $\sigma_r\in S_{2g+2}$ is the permutation
of strands induced by $\hat\tau_r$. Instead of giving the construction
over $S^*$ we construct $B^*$ over the diffeomorphic domain
\[\textstyle
  U\ =\ \Big(B_1(0)\cup\bigcup_{r=1}^\mu B_{2\eps}(s_r)\Big)
  \setminus \bigcup_{r=1}^\mu\overline{B_\eps}(s_r)\,.
\]
We put $A_r(t)=A_r(0)$ for $t<0$, $A_r(t)=A_r(1)$ for $t>1$ and define
\begin{eqnarray*}
  B^*&=&\textstyle \Big(B_1(0)\setminus \bigcup_{r=1}^\mu B_\eps(s_r)\Big)
  \times\{P_1,\ldots,P_{2g+2}\}\\
  &&\cup \bigcup_{r=1,\ldots,\mu \atop k=1,\ldots,2g+2}
  \textstyle \Big\{(s_r+ \rho e^{\pi i (t-\frac{1}{2})+i\arg
s_r},A_r^k(t))\,\Big|\,
  t\in(-\frac{1}{2},\frac{3}{2}),\ \rho\in(\eps,2\eps)\Big\}\,.
\end{eqnarray*}
Now we want to fill in the holes $\overline{B_\eps}(s_r)$. The braid
monodromies $A_r$ around the $s_r$ have been discussed above. We saw
that they are represented either by half-twists of a pair of strands
or
by a double twist of an odd number of strands. In any case, any $A_r$
can be deformed to the monodromy $T$ of one of the standard families
(\ref{standard.curves}) by a family of braids ($s\in[0,1]$)
\[
  A_r[s]:\ [0,1]\ \longrightarrow\ (S_2)^{2g+2}\setminus D\,,\quad
  A_r[0]=A_r, A_r[1]= T\,.
\]
Using this isotopy of braids one can glue in the standard models
(\ref{standard.curves})  to arrive at a possibly singular surface
$\bar B\subset \cz\times S^2$ with
\begin{enumerate}
    \item
    $\bar B\to\cz\setminus \{s_1,\ldots,s_\mu\}$ is an unbranched covering
    \item
    if $q^{-1}(s_r)\subset M$ is irreducible then $B\subset\cz\times
    S^2$ is locally around $s_r$ isomorphic to $B_\mathrm{irr}$
    \item
    if $q^{-1}(s_r)\subset M$ is a union of closed surfaces of genera
    $h$ and $g-h$ ($h\le g-h$), then $\bar B\subset \cz\times S^2$ is
    locally around $s_r$ isomorphic to $B_{h,g-h}$
    \item
    the global monodromy is either trivial or $(x_1\cdot\ldots\cdot
    x_{2g+1})^{2g+2}$.
\end{enumerate}
Here ``locally isomorphic'' refers to the existence of a local
diffeomorphism of
the ambient spaces mapping the branch surfaces to each other.

Next we have to extend $\bar B$ across $\infty$. If the global monodromy
is trivial we use an isotopy of the monodromy over $\{|s|=2\}$ to the
trivial braid $[0,1]\times(P_1,\ldots,P_{2g+2})$ to connect $\bar B\cap
(B_2(0)\times S^2)$ to the product surface $(\cz\setminus
B_4(0))\times\{P_1,\ldots,P_{2g+2}\}$, which readily extends over
$\infty$ in the product bundle $\bar P:= S^2\times S^2\to S^2$.

In the other case of non-trivial monodromy $(x_1\cdot\ldots\cdot
x_{2g+1})$ we may extend over $\infty$ in the non-trivial $S^2$-bundle
$\bar P\to S^2$, which can be realized holomorphically as the first
Hirzebruch surface $\pr(\calo\oplus \calo(1))$: $\bar P$ can be
constructed by patching trivial bundles over half spheres via
\[
  S^1\ \longrightarrow\ \mathrm{Diff}(S^2=\cz\cup\{\infty\})\,,
  \quad e^{2\pi i t}\ \longmapsto\ \Big(z\mapsto e^{2\pi i t}\cdot
  z\Big)\,.
\]
The trivial braid along the equator in the trivialization at $\infty$
clearly extends (as graph of $2g+2$ sections) over the half-sphere
at $\infty$. Applying the patching function to this trivial braid
we obtain however one full twist on the bunch of all strands. We
pointed out above that this braid is nothing but $(x_1\cdot\ldots\cdot
x_{2g+1})^{2g+2}$. This shows that we can extend $B\cap(B_2(0)\times
S^2)$ to a surface in $\bar P$ with no additional branching. We keep
the notation $\bar B$ for the resulting extended surface.

From the local description of the singularities of $\bar B$ we can
desingularize $\bar B$ by two consecutive blow ups of $\bar P$ for each
singular
point. We denote this blow up by $\bar\rho: \tilde P\to \bar P$.
The result is an $S^2$-fibration $\tilde P\to S^2$ with as many
singular fibers of type $\Lambda_2$ as $q:M\to S^2$ has reducible singular
fibers. Let $B\subset \tilde P$ be the strict transform of $\bar B$, and
$B_E$ the union of vertical $(-2)$-spheres that arise in the blow up
process. To produce a two-fold covering $\tilde M\to \tilde P$ with
branch locus $B+B_E$ we only have to show that $B+B_E\in H_2(\tilde P,
\gz)$ is 2-divisible, see Proposition~\ref{twofoldcovers} in the
appendix.
\begin{lemma}
    $B+B_E\in H_2(\tilde P,\gz)$ is 2-divisible.
\end{lemma}
\pf
The class of $B+B_E$ differs by a 2-divisible class from
$\bar\rho^*\bar B$. In fact, if $Q\in\bar B$ is a point of
multiplicity $2h+1$ and $\rho_1:P_1\to\bar P$ is the blow up in $Q$
with exceptional divisor $E\subset P_1$ then
\[
  \rho_1^*\bar B\ =\ B_1+(2h+1)E\,,
\]
where $B_1$ is the strict transform of $\bar B$, cf.\ Fig.~1. Another blow up
$\rho_2:P_2\to P_1$ in the singular point $B\cap E$ with exceptional divisor
$\Gamma'\subset P_2$ and strict transform $B_2$ of $B$ and $\Gamma$
of $E$ leads to
\[
  \rho_2^*\rho_1^*\bar B\ =\ \rho_2^*B_1+(2h+1)\rho_2^*E\
  =\ B_2+(2h+1)\Gamma'+(2h+1)(\Gamma+\Gamma')\,.
\]
This differs from $B_2+\Gamma$ by $2((2h+1)\Gamma'+h\Gamma)$. Doing
this successively at all singular points of $\bar B$ and noticing that
the union of the $\Gamma$ equals $B_E$, we conclude
that it suffices to prove 2-divisibility of $\bar B$.

We now employ the notion of virtual number of critical points
introduced above. Let as before $\bar P=\pr(\calo\oplus \calo(n))$
($n\in\{0,1\}$) and write $\bar B\sim dH+lF$ with $d=2g+2$. We have to show
that
$l\equiv 0\ ({\rm mod 2})$. Proposition~\ref{numbersing} tell us that
\[
  \mu_\mathrm{virt}(\bar B)\ =\ 2(d-1)l+nd(d-1)\,,
\]
so $\mu_\mathrm{virt}(\bar B)-nd(d-1)\equiv 0\ (4)$ iff $l\equiv 0\ (2)$.
We also noted that $\mu_\mathrm{virt}(\bar B)$ is a sum of local
contributions that coincide with the signed length of the word in
$x_1,\ldots,x_{2g+2}$ representing the braid monodromy around
critical values. That this length minus $nd(d-1)$ is divisible by 4
must be related to the fact that our braid monodromy representation comes
from a representation into $\HMC_g$, for there are clearly
surfaces $\bar B\subset\bar P$ homologous to $dH+lF$ with $l$ odd.
We introduce an auxiliary group $\tilde B(S^2,2g+2)$ as fibered product
\[\begin{array}{ccc}
    \tilde B(S^2,2g+2) & \stackrel{/I}{\longrightarrow} & B(S^2,2g+2)  \\
    \diagl{/(x_1\ldots x_{2g+1})^{2g+2}} &  & \diagr{/(x_1\ldots
x_{2g+1})^{2g+2}} \\
    \HMC_g & \stackrel{/I}{\longrightarrow} & \MC(S^2,2g+2)\quad .
\end{array}\]
Explicitely we have to drop the relation $(x_1\ldots
x_{2g+1})^{2g+2}=1$ in $\HMC_g$, or replace the relation
$I=1$ by $I x_i=x_i I$, $I^2=1$ in $B(S^2,2g+2)$. In other words,
\[
  \tilde B(S^2,2g+2)\ =\ \Bigg\langle x_1,\ldots, x_{2g+1}, I
  \Bigg| \begin{array}{c}x_i x_j=x_j x_i\ \ {\scriptstyle |i-j|\ge2}\,;\
           x_i x_{i+1}x_i=x_{i+1}x_i x_{i+1}\\
           I=x_1\ldots x_{2g+1}x_{2g+1}\ldots x_1;\
	   I x_i=x_i I;\ I^2=1\end{array}\Bigg\rangle
\]
As in the case of lifting $\bar\tau_r\in\MC(S^2,2g+2)$ to
$\hat\tau_r\in
B(S^2,2g+2)$ there exists a distinguished lift $\tilde\tau_r$ of $\tau_r$
to $\tilde B(S^2,2g+2)$.  By definition $\tilde\tau_r$ maps to $\hat\tau_r
\in B(S^2,2g+2)$ and
\[
  \tilde\tau_1\cdot\ldots\cdot\tilde\tau_\mu\ =\ \left\{
  \begin{array}{ll}1&,\,n=0\\
      (x_1\ldots x_{2g+1})^{2g+2}&,\,n=1\end{array}\right.
\]
The relations in $\tilde B(S^2,2g+2)$ have the virtue of preserving
the signed length modulo 4. This together with the length
$nd(d-1)$ of the total monodromy  shows that
\[
  \mathrm{length}(\tilde\tau_1\cdot\ldots\cdot\tilde\tau_\mu)
  -nd(d-1)\ \equiv\ 0\quad (4)
\]
as claimed. This concludes the proof of 2-divisibility of $\bar B$,
hence of $B+B_E$.
\qed

The proof of the theorem is now almost finished: By
2-divisibility of $B+B_E$ there exists a two-fold cover
$\tilde\kappa:\tilde M \to\tilde P$ branched along $B+B_E$. Every
connected component $\Gamma$ of $B_E$ is a sphere of
self-intersection $-2$. Hence $\kappa^{-1}(\Gamma)$ is a
$(-1)$-sphere and can be contracted. Let $\sigma:\tilde M\to M$ be
the blow down of all these $(-1)$-spheres. A local study of this
process shows that every contracted component results in an ordinary
critical point of the induced fibration $q':M\to S^2$, which is thus
a SLF. By Lemma~\ref{compatibility.of.mon} and the construction, the
monodromy of $q'$ coincides with the one coming from $q:M\to S^2$.
The result of Kas implies that $q$ and $q'$ are isomorphic SLF's.

Finally, it is well known (and can be easily verified) that the
contraction of a $(-2)$-sphere $\Gamma$ leads to the orbifold
singularity
\[
  \cz^2/(\gz/2\gz)\ \simeq\ \{(z,w,y)\in\cz^3\,|\,z^2+w^2+y^2=0\}\,.
\]
This shows that $q:M\to S^2$ factorizes over the orbifold $P$ as
described in statement 4 of the proposition.

The converse assertion, that a SLF of this form is hyperelliptic, is
trivial.
\qed

\begin{remark}\rm
    The theorem shows that hyperelliptic SLF's $M\to S^2$ admit a
    hyperelliptic involution $\tau:M\to M$ with fix locus the
    embedded surface $\sigma(\tilde\kappa^{-1}(B))$ union the nodes of
    reducible singular fibers ($=\sigma(\tilde\kappa^{-1}(B_E))$\,).
    This may fail if one admits multiple singular fibers.
\qed
\end{remark}
A SLF $q:M\to S^2$ may of course admit several {\em hyperelliptic
structures}, that is, choices of two-fold covers $\Sigma=q^{-1}(s)
\to S^2$ commuting with the monodromy up to a diffeomorphism
of $S^2$ and up to isotopy. Essentially different choices lead to
non-isomorphic hyperelliptic involutions $\tau: M\to M$, hence to
non-isomorphic pairs $\tilde B\subset P$. An example
is provided by elliptic Lefschetz fibrations: Moishezon and Livn\'e
have shown that elliptic Lefschetz fibrations are determined by their
number of singular fibers $\mu$: It is always
$\mu\equiv 0(12)$ and the monodromy can be put into the form
$(x_1, x_2)^{\mu/2}$ \cite{moishezon}. But the case $\mu=12$
arises either with $B\sim 4H+2F$ in $\pr^1\times\pr^1$ or with
$B\sim 4H$ in $\pr(\calo\oplus\calo(1))$ (cf.\ Theorem~\ref{invts}).

The theorem motivates the following definition.
\begin{definition}\rm\label{branch.surf}
  Let $\tilde p:\tilde P\to S^2$ be an $S^2$-fibration with $t$ singular
  fibers of type $\Lambda_2$ and $\bar P\to S^2$ the associated $S^2$-bundle,
  $\bar \rho:\tilde P\to \bar P$. An embedded surface $B\subset \tilde P$ will
  be called {\em branch surface} (of type $(g,\mu,t)$\,) if
  \begin{itemize}
      \item
      the composition $B\to\tilde P\to S^2$ is a branched cover of
      degree $2g+2$ with $\mu-t$ ordinary critical points
      (that is, of the form $z\mapsto z^2$ in appropriate complex
      coordinates)
      \item
      let $\tilde p^{-1}(s)=F\cup\Gamma\cup\Gamma'$ with $F^2=-1$,
      $\Gamma^2=-2$, ${\Gamma'}^2=-1$ be a singular fiber. Then
      $B\pitchfork\tilde p^{-1}(s)$, $B\cap\Gamma=\emptyset$, $B\cdot
      F\equiv 1\ (2)$
      \item
      $\bar\rho_*[B]\in H_2(\bar P,\gz)$ is 2-divisible
      \qed
  \end{itemize}
\end{definition}
The theorem can then more concisely be phrased as
\begin{corollary}
    For $g\ge2$ there is a one-to-one correspondence between isomorphism
    classes of SLF's with hyperelliptic structure $(q:M\to S^2,\kappa:
    q^{-1}(s) \to S^2)$, of genus $g$ with $\mu$ singular fibers, $t$
    of which are reducible, and isomorphism classes of pairs $(\tilde
    P,B)$ with $\tilde P\to S^2$ an $S^2$-fibration with $t$ fibers of
    type $\Lambda_2$ and $B\subset \tilde P$ a branch surface of type
    $(g,\mu,t)$. Moreover, for $g=2$ the hyperelliptic structure
    is determined by the fibration up to isomorphy.
\end{corollary}
\pf
Starting from a branch surface $B\subset\tilde P$ the two-fold cover
branched along $B\cup B_E$ is a hyperelliptic SLF $q:M\to S^2$. The proof
of the theorem now produces a possibly singular branch surface $\bar B'$
in the $S^2$-bundle $\bar P$ associated to $\tilde P$ (the isomorphism
class of which was determined by the monodromy at infinity).
There was one non-uniqueness involved in the construction of
$\bar B'$. To wit, if $g=2h$ is even then for any reducible singular
fiber of type $B_{h,h}$ we can choose either $(x_1\cdot\ldots
\cdot x_h)^{4h+2}$ or $(x_{h+2}\cdot\ldots
\cdot x_g)^{4h+2}$ as word representing the corresponding
Dehn-twist. This amounts to select one of the two connected components
of $S^2\setminus\kappa(C)$ if $C\subset \Sigma$ is the $I$-invariant
embedded circle supporting the Dehn-twist. However, both choices lead to
the same branch surface after blow up, where this non-uniqueness
materializes in a choice of one of the two $(-1)$-curves to be
contracted first. This shows that the desingularization $B'\subset
\tilde P$ of $\bar B'$ is isomorphic to $B\subset P$.

For $g=2$ any cover $\kappa:\Sigma\to S^2$ is compatible with
the monodromy and hence gives rise to a hyperelliptic structure.
But any two such covers $\kappa, \kappa'$ with the same branch set
in $S^2$ differ only by a diffeomorphism of $\Sigma$. But the
mapping class group $\MC_2=\HMC_2$ is an extension of the mapping
class group of $(S^2, \{P_1,\ldots,P_6\})$ by the hyperelliptic
involution. This shows that for $g=2$ there is indeed only one
hyperelliptic structure on $\Sigma$ up to isomorphy.
\qed
\vspace{1ex}

Note that in contrast to $B\subset\tilde P$ the possibly singular branch
surface $\bar B$ in the $S^2$-bundle $\bar P$ associated to $\tilde P$ is
not unique.
\begin{remark}\rm
  The construction of the two-fold cover from the possibly singular surface
  $\bar B\subset \bar P$ suggests the question why ordinary singular points
  do not occur in $\bar B$. For multiplicity $h$ these have typical braid
  monodromy $(x_1\cdot\ldots\cdot x_{h-1})^h$ (one full twist of the
  first $h$ strands) and they require only one blow up for
  desingularization, so they seem to be more basic than infinitely
  close singular points. But then $\tilde M\to \tilde P$ is unbranched
  over the critical point of $\tilde P\to S^2$ (the intersection of
  the exceptional divisor of $\tilde P\to\bar P$ with the strict
  transform of the fiber of $\bar P\to S^2$). Thus $M\to S^2$
  acquires a singular fiber with two critical points that connect two
  smooth components of genera $\frac{h}{2}-1$ and $g-\frac{h}{2}$
  ($h\equiv0$ is forced by 2-divisibility of $B$). The monodromy of
  $M$ is therefore the product of two commuting Dehn-twists of
  non-disconnecting type.

  It is also interesting to note the result of a slight (locally
  holomorphic) perturbation of the projection $q:M\to S^2$. It can
  be shown that for ordinary singular points in $\bar B$
  this is equivalent to smoothing $\bar B$
  inside $\bar P$. This is an example of
  Brieskorn's simultaneous resolution \cite{brieskorn}, which suggested
  the notion of {\em inessential singularities} of branch loci in
  algebraic geometry \cite{horikawa}\cite{persson}. In the case of an ordinary
  singular point of multiplicity $h$ we obtain $h(h-1)$ critical points
  of type $B_\mathrm{irr}$, so the ``stabilization'' of the fibration
  belonging to an ordinary singularity of multiplicity $h$ is a
  fibration with $h^2-h$ irreducible singular fibers.
\end{remark}


\section{Study of the branch locus}
We now want to collect more information on the (non-vertical
components of) branch loci $B\subset\tilde P$, or, equivalently $\bar B
\subset\bar P$. We begin by showing that $B$ is symplectic. Without
more effort this can be proved for SLF's factorizing (possibly after blow
up) over an arbitrary branched covering and the blow-up
of a symplectic fiber bundle.
\begin{theorem}\label{B_sympl}
    We assume that $q:M\to S^2$ is a SLF and that the induced fibration
    of a blow up $\sigma:\tilde M\to M$ factorizes over a branched
    covering $\tilde\kappa$ (any degree, any branch type)
    \[
      \tilde M\stackrel{\tilde\kappa}{\longrightarrow} N
      \stackrel{p}{\longrightarrow} S^2
    \]
    with $p: N\to S^2$ the blow up of a symplectic fiber bundle. Let
    $B\subset N$ be the union of the non-vertical components of the
    branch locus and assume $B\cap\mathrm{Crit}(p)=\emptyset$.
    Let $\omega_N$ be a symplectic form on $N$ such that the general
    fibers of $p$ are symplectic and denote by $\omega_{S^2}$ a
    symplectic form on $S^2$.

    Then $B$ is symplectic with respect to $\omega_N+k\cdot
    p^*\omega_{S^2}$ for sufficiently large $k$.
\end{theorem}
\pf
The composition $p|_B:B\hookrightarrow N\to \pr^1$ is a branched
covering. We endow $B$ with the unique complex structure making
this covering holomorphic. Then $p^*\omega_{\pr^1}|_B$ is a K\"ahler
form on $B$ away from the critical set, so $\omega_N+k\cdot
p^*\omega_{\pr^1}|_B$ will be non-degenerate away from arbitrarily
small neighbourhoods of the critical set. It remains to show that for
each critical point $P\in B$, the restriction $\omega_N$ to $T_P B$
is non-degenerate and agrees with the complex orientation.

By holomorphicity of $p|_B$ there are local holomorphic coordinates
$w$ for $B$ near $P$ and $t$ for $\pr^1$ with $(p|_B)^*t=w^m$ for
some $m\ge2$. Extend $w$ to a smooth function on a neighbourhood of
$P$ in $N$. Since $p_*T_PB=0$, the tuple $(z:=p^*t,w)$ is a complex
coordinate system on $N$, and $\omega_N|_{T_PB}$ is non-degenerate.
Now if $Q\in \kappa^{-1}(P)$ is a branch point of $\tilde\kappa$ of order
$b\ge 2$, then $M$ can be described near $Q$ by
\[
  \{(z,w,u)\in\cz^3\,|\, u^b=z-w^m\}\,.
\]
So $(u,w)$ are complex coordinates on $M$ near $Q$ and $p:(u,w)\mapsto
u^b+w^m$. Since by hypothesis $q:M\to\pr^1$ is oriented and Lefschetz,
we see that necessarily $b=m=2$ and $(u,w)$ is compatible with the
orientation of $M$. In turn $(z,w)$ agrees with the orientation of
$N$ (and of $\pr^1$). The equality of complex line bundles
\[\textstyle
  T_PB\ =\ \ker p_*\ =\ \det_\cz T_P N\otimes_\cz p^* K_{\pr^1}
\]
together with the fact that $p:N\to\pr^1$ is symplectic with respect
to $\omega_N$ then shows that the orientations of $T_PB$ given by
$\omega_N|_{T_PB}$ on one hand and by $w$ on the other hand coincide.
\qed

\begin{remark}\label{symplectic.rem}\rm
The authors believe that symplecticity of the branch locus is an essential
observation, because it reduces the possible isomorphy classes of branch
loci drastically. There are in fact good chances that in $S^2$-bundles over
$S^2$, and in $\pr^2$, symplectic submanifolds are classified even up to
isotopy by their homology, just as in the holomorphic situation. In view of
Theorem~\ref{red.to.branch} this would in particular imply that every
hyperelliptic SLF with only irreducible singular fibers is (isomorphic  to) a
holomorphic (hyperelliptic, Lefschetz) fibration. At least for surfaces inside
$S^2$-bundles over $S^2$ of bidegree $(2,d)$ or $(3,d)$ this can indeed be
shown by studying the braid monodromy. Without symplecticity there
would be no hope for a simple classification, as one can always paste
non-trivial 2-knots locally without changing the homology class. Our
approach to the isotopy problem of symplectic submanifolds also suggests
that holomorphicity of hyperelliptic fibrations is true provided the  number
$\mu$ of singular fibers is large compared to the number $t$  of reducible
singular fibers (a sufficient condition might be $\mu/t>  18$), so reducible
singular fibers might be interpreted as  obstructions to holomorphicity. In
fact, known examples of  non-holomorphic, hyperelliptic SLF's have small
ratio $\mu/t$, e.g.\ $\mu/t=4$ for the recent example in
\cite{ozbagci}.
\qed
\end{remark}
\vspace{2ex}

\noindent
We can also prove that unless $M\to S^2$ is trivial, $B$ has at most
two connected components, and if $B$ is disconnected one of them
can be taken as negative section $S_\infty$ of
$\pr(\calo\oplus\calo(m))$ with $m$ even.
\begin{proposition}\label{connected}
  Let $B\subset \tilde P$ be a branch surface
  (Definition~\ref{branch.surf}). Then either
  $B$ is connected or the reduction $(\bar P,\bar B)$ is diffeomorphic
  (as $S^2$-bundle over $S^2=\pr^1$) to
  \begin{itemize}
  \item
    $(\pr(\calo\oplus\calo(2n)),S_\infty\cup \bar B')$ with $n>0$ and
    $\bar B'\sim (2g+1)H$ connected
  \item
    $(\pr^1\times\pr^1, H_1\cup\ldots\cup H_{2g+2})$ where $H_i$ are disjoint
    sections of  $P$
  \end{itemize}
\end{proposition}
\pf
By the local study of the braid monodromy it suffices to investigate
connectedness of a smoothing of $\bar B\subset\bar P$. In other words,
we may assume that $P=\bar P$ and $B=\bar B$ is already smooth.
Let us first consider the case $P\simeq \pr^1\times\pr^1$. If $B$ has more than
two components, say $B_i\sim a_i H+b_i F$, $i=1,2,3$ (necessarily $a_i>0$),
then
$b_i/a_i=0$ for the pairwise intersection is zero. Now Lemma~\ref{numbersing}
applied to $B_i$ shows that $B_i\rightarrow \pr^1$ is unbranched, hence $a_i=1$
for all $i$ and we get case (2).

If $B$ has two components we may write them as $S'\sim aH+bF$, $B'\sim
cH+dF$ with
$a,c>0$, $a+c=2g+2$, not both $b$ and $d$ equal to 0, and ordered in such a
way that $(S')^2= 2ab \le (B')^2=2cd$. Now $0=S'\cdot B'=ad+bc$, so $b\le0$
and $d\ge0$. Moreover, since $S'$ is symplectic, the adjunction formula shows
\[
  2g(S')-2\ =\ (-2H-2F)(aH+bF)+(S')^2\ =\ -2b-2a+2ab\,,
\]
hence $g(S')=-b(1-a)-a+1$. But $a\ge 1$ and $b\le0$ so $g(S')\le0$, which
implies $g(S')=0$, $a=1$, $d=-bc=-b(2g+1)$. Lemma~\ref{disconncase}
below thus shows that we are dealing with case (2) with $n=-b$.

The other possibility is $P$ diffeomorphic (as fiber bundle)
to $\pr(\calo\oplus\calo(1))$. If $B$
were not connected write $B=B_1\cup B_2$ with $B_1^2\le B_2^2$, $B_1=a
S_\infty+bF$,
$B_2=c H+dF$, $a,c>0$, $a+c=2g+2$. The adjunction formula for $B_1$ yields
\[
  2g(B_1)-2\ =\ (-2S_\infty-3F)(aS_\infty+bF)+(aS_\infty+bF)^2\ =\
  -(a-b)^2-a\ \le\ -a\ <\ 0\,,
\]
so $g(B_1)=0$, $a=1$ and $b\in\{0,1\}$. If $b=1$ the equation $0=B_1\cdot
B_2=ad+bc$
implies $d=-c$ and then $B_2^2=-c^2$, $B_1^2=1$ contradicts $B_1^2\le B_2^2$.
The remaining case is $B_1=S_\infty$, $B_2=(2g+1) H$. But then $B_1+B_2=(H-F)+
(2g+1)H=(2g+2)H-F$ is not a 2-divisible class and can thus not occur as branch
locus of a two-fold covering.
\qed
\vspace{1ex}

\noindent
We owe one lemma.
\begin{lemma}\label{disconncase}
  Suppose that $P\to \pr^1$ is a ruled surface and $S\subset P$ is a
  section with $S^2=-n<0$. Then
  there exists an isomorphism of fiber spaces $P\simeq
  \pr(\calo\oplus\calo(n))$ carrying
  $S$ to the negative section $S_\infty$.
\end{lemma}
\pf
Removing $S$ we obtain a fiber bundle $P\setminus S\rightarrow \pr^1$
with structure group ${\rm Aff}(\cz)=\cz\rtimes \cz^*$, the complex
affine transformations of $\cz$. Conversely, fiber bundles with fiber
$(\cz,{\rm
Aff}(\cz))$ correspond exactly to pairs $(P,S)$ of $\pr^1$-bundles with
section. Now
write $P\setminus S$ as two copies of $\Delta\times\cz$ glued over an annulus
$A\subset \Delta$ via a map $\theta: A\rightarrow {\rm Aff}(\cz)$. But
$\cz^*\subset{\rm Aff}(\cz)$ is a deformation retract and by degree theory
any map
$A\rightarrow\cz^*$ is homotopic to a holomorphic map, in fact to $\theta_k:
z\rightarrow z^k$ for some $k\in\gz$. So after a change of trivialization
we may
assume $\theta=\theta_k$. But gluing two copies of $A\subset \Delta$ via
$\theta_k$  we obtain $\pr(\calo\oplus\calo(k))\setminus S_\infty$ and this
extends to an isomorphism $P\simeq
\pr(\calo\oplus\calo(k))$ mapping $S$ to $S_\infty$. The number $k$ is
determined by $k=-(S_\infty)^2=-S^2=n$.
\qed
\begin{remark}\rm
In two recent papers Fuller claimed that any $SLF$ can be factorized into a
simple 3-fold branched cover of an $S^2$-fibration over $S^2$ followed by
the bundle projection \cite{fuller1}, \cite{fuller2}. Unfortunately, our
observation on symplecticity of the branch locus shows that this is
impossible: The covering being simple forces the branch locus to
decompose into two parts, each of which having covering degree at least
two over $S^2$. As we saw in the proof of Proposition~\ref{connected} this
is impossible for homological reasons unless the branch locus is trivial. The
problem in op.cit.\ is that the branch locus, which is first constructed in
$\Delta\times S^2$, can not be extended over the fiber over $\infty\in
S^2$.
\qed
\end{remark}


\section{Topological invariants of hyperelliptic SLF's}
We will see in this section that the description as relative minimalization
of a branched cover allows the computation of topological invariants of
hyperelliptic SLF's in  terms of the genus $g$ and the types and numbers of
singular fibers.

Let as before $q:M\to S^2$ be a hyperelliptic SLF of genus $g$ with
$\mu$ singular fibers, $t$ of which are reducible of types
$h_1,\ldots,h_t$ (the minimum of the genera of the two irreducible
components). Let $\sigma:\tilde M\to M$ be the blow up of the $t$
nodes of reducible singular fibers, $\tilde \kappa:\tilde M\to
N:=\tilde P$ the two-fold cover with branch locus $B_N=B\cup B_E$,
and $\bar\rho: N\to \bar P$ the birational map
to an $S^2$-bundle $\bar P=\pr(\calo\oplus\calo(n))$. Some more notations
for later reference: $E_1,\ldots,E_t\subset \tilde M$ are the
exceptional fibers of $\sigma$, $\Gamma_1\cup\Gamma'_1,\ldots,
\Gamma_t\cup\Gamma'_t$ are the exceptional fibers of $\bar\rho$ with
$\Gamma_i^2=-2$, ${\Gamma'_i}^2=-1$, $\bar B=\bar\rho(B)$, $\tilde
B=\kappa^{-1}(B_N)$.

For the following it will be convenient to work in the almost complex
category as follows: We know that $B\cup B_E\subset N$ is symplectic
(with respect to $\omega'_N:=\omega_N+k\cdot p^*\omega_{\pr^1}$
for $k\gg 0$: Theorem~\ref{B_sympl}), so there exists an
$\omega'_N$-tamed almost complex structure $J$ on $N$ making $B\cup
B_E$ a $J$-holomorphic curve. We may also assume that $J$ is
compatible with the fiber structure, so $N\to\pr^1$ is a holomorphic
map of almost complex manifolds (near the critical values we can
choose $B$ holomorphic with respect to an integrable complex
structure). Clearly, $J$ lifts to an almost complex structure $\tilde
J$ on $\tilde M$ such that $\kappa:\tilde M\to N$ is $(\tilde
J,J)$-holomorphic. By the standard construction of a symplectic structure
$\tilde\omega$ on branched covers (unique up to deformation
equivalence), $\tilde J$ is $\tilde \omega$-tame for an appropriately
chosen $\tilde \omega$ on $\tilde M$.

The advantage of introducing (tamed) almost complex structures is that
we may now employ almost complex analogues of well-known formulas
from the theory of complex surfaces. To give them the usual shape, we
introduce the following notation for topological invariants of an
almost complex four-manifold $(M,J)$:
\[
  K_M\ :=\ -c_1(M)\,,\quad\chi(M)\ :=\ \frac{c_1^2(M)+e(M)}{12}\,,
\]
where $c_i(M)$ denotes the $i$-th Chern class of $(T_M,J)$ and
$e(M)=c_2(M)$ is the Euler class. We identify $H^4(M,\gz)=\gz$
and $H_2(M,\gz)$ with $H^2(M,\gz)$ by Poincar\'e-duality, as is
customary in the theory of complex surfaces.
With this understood we have in our situation
\begin{eqnarray}\label{adjunction}\begin{array}{rcl}
    K_{\tilde M}&=&\sigma^* K_M+\sum_i E_i\ ,
    \ \mathrm{so}\quad K_M^2\ =\ K_{\tilde M}^2+t\\[1ex]
    K_{\tilde M}&=&\kappa^* K_N+\tilde B\ =\
    \kappa^*(K_N+\frac{1}{2}B_N)\ ,
    \quad\mathrm{so}\ K_{\tilde M}^2=2(K_N+\frac{1}{2}B_N)^2 \\[1ex]
    \bar\rho^*\bar B&=&B+\sum_i(2h_i+1)(\Gamma_i+2\Gamma'_i)\\[1ex]
    &=& B+B_E+\sum_i(2h_i\Gamma_i+(2h_i+1)2\Gamma'_i)\\[1ex]
    K_N&=&\bar\rho^*K_{\bar P}+\sum_i(\Gamma_i+2\Gamma'_i)
\end{array}\end{eqnarray}
We are now ready to compute $e$ and $c_1^2$.
\begin{theorem}\label{invts}
    Let $q:M\to S^2$ be the hyperelliptic SLF of genus $g$ with $\mu$
    singular fibers, $t$ of which are reducible, of types $h_1,\ldots,h_t$,
    associated to a branch locus $\bar B\subset \bar P=\pr(\calo\oplus
    \calo(n))$ homologous to $(2g+2)H+2kF$. Then
    \begin{eqnarray*}
      \mu&=&(2g+1)(4k+n(2g+2))-{\textstyle\sum_i}(8h_i^2+4h_i-1)\\[1ex]
      e(M)&=&4-4g+\mu\\[1ex]
      c_1^2(M)&=&(g-1)(4k+n(2g+2)-8)-
      {\textstyle\sum_i}(2h_i-1)^2\\[1ex]
      \tau(M)&=&-(g+1)(4k+n(2g+2))+{\textstyle\sum_i}
      (4h_i^2+4h_i-1)\\[1ex]
      \chi(M)&=&\frac{g}{4}(4k+n(2g+2))-g+1-{\textstyle\sum_i} h_i^2
    \end{eqnarray*}
\end{theorem}
\pf
By the discussion following the definition of $\mu_\mathrm{virt}$ we
have
\[
  \mu_\mathrm{virt}(\bar B)\ =\ (\mu-t)+\sum_{i=1}^t 4h_i(2h_i+1)\,,
\]
but also $\mu_\mathrm{virt}(\bar B)=2(2g+1)2k+n(2g+2)(2g+1)$
by Lemma~\ref{numbersing}, which yields the formula for $\mu$.

The formula for the Euler number can be proven by taking $U\subset S^2$ to
be an $\eps$-neigh\-bour\-hood of the critical locus and noting
\[
  e(M)\ =\ e(q^{-1}(S\setminus U))+e(q^{-1}(\bar U))\,.
\]
Over $S\setminus U$ the fibration is locally trivial, so the first term gives
$e(S\setminus U)\cdot e(F)$, which for $S=\pr^1$ equals
$(2-\mu)(2-2g)$. On the other hand $q^{-1}(\bar U)$ is a deformation
retract of the set of singular fibers, and each of these having only
one node (that we think of as contracted circle on a general fiber)
contributes $3-2g$. Putting things together we obtain the stated
formula.

As for $c_1^2(M)=K_M^2$, using (\ref{adjunction}) the computation
is straightforward: The formulas for $K_N$ and $B_N$ yield
\[
  K_N+\frac{1}{2}B_N\ =\ \bar\rho^*(K_{\bar P}+\frac{1}{2}\bar B)-
  \sum_i\Big((2h_i-1)\Gamma'_i+(h_i-1)\Gamma_i\Big)\,,
\]
so we get
\begin{eqnarray*}
 \lefteqn{\Big(K_N+\frac{1}{2}B_N\Big)^2}\hspace{1cm}\\
 &=&\Big(-2H+(n-2)F+(g+1)H+kF\Big)^2+
    {\textstyle\sum_i}\Big((2h_i-1)\Gamma'_i+(h_i-1)\Gamma_i\Big)^2\\
  &=&(g-1)(2k+n(g+1)-4)+{\textstyle\sum_i}(-2h_i^2+2h_i-1)\,.
\end{eqnarray*}
But also
\[
  K_M^2\ =\ K_{\tilde M}^2+t\ =\ 2(K_N+\frac{1}{2}B_N)^2+t\,,
\]
from which the result follows. Index $\tau$ and ``holomorphic Euler
characteristic'' $\chi$ can  be computed in terms of $c_1^2$ and
$e$ as noted before.
\qed

Observe that $n$ and $k$ occur only in the combination $4k+n(2g+2)$
and can therefore be eliminated from the formula for $c_1^2$ in terms
of $\mu$, $h_i$ and $g$.
\vspace{2ex}

The fundamental group of a SLF is a quotient of the fundamental group
of a general fiber and can generally be arbitrary \cite{katzarkovetal}.
However, under the absence of reducible singular fibers one
can show simply-connectedness of hyperelliptic SLF's.
\begin{proposition}\label{fund.group}
    Non-trivial, hyperelliptic SLF's without reducible fibers are
    simply-connected.
\end{proposition}
\pf
The handlebody decomposition of $M$ associated to a DLF $q:M\to
S^2$ \cite{kas} implies the following well-known description of the
fundamental group of $M$: Letting $\delta_1,\ldots,\delta_\mu:
S^1\hookrightarrow \Sigma$ be the vanishing cycles, it holds
\[
  \pi_1(M)\ =\ \pi_1(\Sigma)/(\delta_1,\ldots,\delta_\mu)\,.
\]
A more convenient interpretation is as the fundamental group of the
topological space $\Sigma/(\delta_1,\ldots,\delta_\mu)$ obtained
from $\Sigma$ by contracting the vanishing cycles (Seifert-Van Kampen).

Now if $q:M\to S^2$ is hyperelliptic and without reducible fibers,
let $M\to P=\pr(\calo\oplus\calo(n))\to\pr^1$ with branch locus
$B\subset P$ be the factorization over a two-fold cover as in
Theorem~\ref{red.to.branch}. We can choose the $\delta_i$ as lifts
under $\kappa:\Sigma\stackrel{2:1}{\to}S^2$ of segments
$\bar\delta_i: [0,1]\to S^2$ joining a pair of critical points of
$\kappa$, cf.\ the discussion following Lemma~\ref{comp.van.cycles}.
Then $\Sigma/(\delta_1,\ldots,\delta_\mu)$ is a
two-fold cover of $S^2/(\bar\delta_1,\ldots,\bar\delta_\mu)$. The
latter space is a wedge of spheres. Moreover, since the branch locus
$B$ has at most two connected components
(Proposition~\ref{connected}), the map
\[
  \kappa':\ \Sigma/(\delta_1,\ldots,\delta_\mu)\ \longrightarrow
  \ S^2/(\bar\delta_1,\ldots,\bar\delta_\mu)
\]
is branched only in the vertex of the wedge and in at most one more
point. But any connected, two-fold cover of $S^2$ branched in at most
two points is simply connected. Thus $\pi_1(M)=\pi_1 (\Sigma)/
(\delta_1,\ldots,\delta_\mu)=0$ as claimed.
\qed


\section{Symplectic Noether-Horikawa surfaces}
We conjectured in Remark~\ref{symplectic.rem} that hyperelliptic
SLF's with only irreducible singular fibers are holomorphic. In this
chapter we discuss some more properties of this class of symplectic
four-manifolds. We mostly restrict ourselves to the case $g=2$, the
reason being that in the holomorphic situation this class of complex
surfaces (admitting Lefschetz fibrations of genus 2 without reducible
fibers up to deformation) are distinguished in the theory of
complex surfaces by their topological invariants. Moreover, there is
a simple classification up to algebraic deformations: It consists of
3 infinite series, and each single class is sweeped out by
one algebraically irreducible family. As in the topological situation above,
they are two-fold covers of rational surfaces. The precise result is:
\begin{theorem}\label{horikawa.surf}
    \textnormal{(\cite{horikawa},\cite[p.51]{chakiris})}
    Let $q:M\to\pr^1$ be a non-trivial genus 2 fibration
    (holomorphic with smooth total space, but not necessarily Lefschetz)
    without reducible fibers. Then $q$ can be (holomorphically)
    deformed to a Lefschetz fibration
    \[
      q':\ M\ \stackrel{\kappa}{\longrightarrow}\
      \pr(\calo\oplus\calo(n))\ \stackrel{p}{\longrightarrow}\ \pr^1\,,
    \]
    where $\kappa$ is a two-fold cover branched along a smooth curve
    $B\subset\pr(\calo\oplus\calo(n))$ and
    such that one of the following holds
    \begin{tabbing}
    \hspace*{2ex}\=\textnormal{(I$_k$)}\hspace{3ex}\=
      $n=0$, $B\sim 6H+2kF$ is connected, $k>0$\\[1.5ex]
    \>\textnormal{(II$_k$)}\>
      $n=1$, $B\sim 6H+2kF$ is connected, $k\ge0$\\[1.5ex]
    \>\textnormal{(III$_k$)}\>
      $n=2k>0$, $B=S_\infty\cup B'$ with
      $B'\sim 5H$ connected
    \end{tabbing}
    Moreover, each of these classes comprises exactly one deformation
    type with irreducible parameter space.
\qed
\end{theorem}
Using specialization to reducible branch curves (e.g.\ $6$ sections,
each linearly equivalent to $H$, plus $2k$ fibers in case I$_k$), it
is not hard to compute the braid monodromies. The result in terms of
our standard generators $x_1,\ldots,x_5$ of $B(S^2,6)$ is
\cite[p.115]{chakiris}
\begin{tabbing}
    \hspace*{2ex}\=\textnormal{(I$_k$)}\hspace{3ex}\=
      ($\mu=20k$)\hspace{8ex}\=$(\tau_1,\ldots,\tau_\mu)=$\=
      $(x_1,x_2,x_3,x_4,x_5,x_5,x_4,x_3,x_2,x_1)^{2k}$\\[1.5ex]
    \>\textnormal{(II$_k$)}\>
      ($\mu=20k+30$)\>$(\tau_1,\ldots,\tau_\mu)=$\>
      $(x_1,x_2,x_3,x_4,x_5,x_5,x_4,x_3,x_2,x_1)^{2k}$\\
      \>\>\>\>$(x_1,x_2,x_3,x_4,x_5)^6$\\[1.5ex]
    \>\textnormal{(III$_k$)}\>
      ($\mu=40k$)\>$(\tau_1,\ldots,\tau_\mu)=$\>
      $(x_1,x_2,x_3,x_4)^{10k}$
\end{tabbing}
where $\tau_r$ is the braid monodromy along $\gamma_r$
($\gamma_1,\ldots,\gamma_\mu$ a standard generating set
for $\pi_1(S^2\setminus\{t_1,\ldots,t_\mu\})$ as above).
The explicit form shows that all such surfaces arise as fiber
connected sums of the three basic fibrations I$_1$, II$_0$ and
III$_2$ (the fiber connected sum is unique except in the last case
because the monodromy generates the whole mapping class group).

From the monodromy representation one derives also
easily that $M$ is simply connected (cf.\ below). Since we
know signature $\tau$ and Euler characteristic $e$,
it remains to investigate the parity of the intersection form to
classify these complex surfaces up to homeomorphy. It turns out that
$M$ is spin (has even intersection form) only for type III$_k$ with
$k\equiv 0\ (2)$ \cite{perssonetal}.

From a complex surface point of view these surfaces are minimal
except in cases II$_0$ and III$_2$ and of general type provided $\mu>40$
(i.e.\ excluding types I$_1$, I$_2$, II$_0$, III$_2$). Moreover, they
fulfill
\[
  c_1^2(M)\ \equiv\ 0\ \ (2)\,,\quad
  \chi(M)\ =\ \frac{1}{2}c_1^2(M)+3\,,
\]
so they lie on the Noether line. It is also worthwile to discuss the
exceptional cases as they are the basic building blocks:
\begin{itemize}
    \item
    I$_1$: Projecting onto the second factor inside $\pr^1\times\pr^1$
    we obtain a $\pr^1$-fibration over $\pr^1$ with $12$ singular
    fibers, so this surface is rational.
    \item
    I$_2$: Similarly we obtain an elliptic fibration over $\pr^1$
    with $36$ singular fibers; this is (hence) a fiber conected sum
    of three copies of the basic elliptic fibration of Livn\'e/Moishezon
    mentioned above.
    \item
    II$_0$: Since $S_\infty$ is disjoint from the branch locus, $M$
    contains two $(-1)$-curves. Contracting these we obtain a
    two-fold cover of $\pr^2$ branched along a sextic, which is
    a K3-surface.
    \item
    III$_2$: After contraction of the $(-1)$-curve over $S_\infty$
    one obtains a surface with $c_1^2>0$ which is (thus) of
    general type.
\end{itemize}
Conversely, any minimal complex surface with these invariants
is of this form, or a two-fold cover  of $\pr^2$ branched along
a curve of degree $8$ or $10$ \cite{horikawa}. There is a similar
classification for $c_1^2$ odd. Surfaces on the Noether-line are
therefore sometimes called {\em Noether-Horikawa surfaces}.
\vspace{2ex}

Let us now come back to symplectic geometry.
We suggest the term {\em symplectic Noether-Horikawa surfaces}
for minimal symplectic surfaces on the Noether line that are (symplectically)
birational to a (symplectic) two-fold cover of a rational surface
($\pr^2$ or a  blow-up of $\pr(\calo\oplus\calo(n))$\,). If $M$ is a
cover of $\pr^2$ branched along a symplectic surface of degree $2d$,
a computation as in the previous chapter shows
\[
  c_1^2(M)\ =\ 2d^2-12d+18\,,\quad
  e(M)\ =\ 4d^2-6d+6\,,\quad
  \chi(M)\ =\ \frac{d^2-3d+4}{2}\,,
\]
so $M$ is on the Noether line iff $d=4$ ($c_1^2=2$, $\chi=4$) or
$d=5$ ($c_1^2=8$, $\chi=7$). Otherwise one can use the technique of
pseudo-holomorphic curves to show that a branched cover of a rational
ruled surface is a hyperelliptic SLF without reducible singular fibers.
Conversely, it follows immediately from Theorem~\ref{invts} that a
hyperelliptic SLF is a symplectic Noether-Horikawa surface
iff all singular fibers are irreducible.
\begin{proposition}
    Let $p:P\to \pr^1$ be a rational ruled surface and $\kappa:\tilde
    M\to P$ a two-fold cover branched along a symplectic
    submanifold $B\subset P$. Then $p$ can be deformed to an
    $S^2$-fibration $p':P\to\pr^1$ such that $q=p'\circ \kappa:\tilde
    M\to \pr^1$ is a (hyperelliptic) SLF without reducible singular
    fibers.
\end{proposition}
\pf
(sketch)\ \ Choose a tamed almost complex structure $J$ on $P$ making
$B$ a $J$-holomorphic curve. Let $S$, $F\subset P$ be a section and a
fiber. One can show that for any tamed $J$, the moduli space of
(not necessarily irreducible, but connected) $J$-holomorphic curves
$C$ homologous to $F$ is isomorphic to $S$, simply by sending $C$ to the
unique point of intersection with $S$. Since these $J$-holomorphic
curves are pairwise disjoint for homological reasons, and rational by
the adjunction formula, we obtain another $S^2$-fibration
$p':P\to S\simeq \pr^1$. It can be deformed to $p$ by connecting $J$
to the integrable complex structure within the space of tamed almost
complex structures. For all this we refer to \cite[\S4]{mcduff} for details.
Moreover, by choosing the path $\{J_t\}_t$ generic, one can arrange the
family of fibers to be transversal to the family of $J_t$-holomorphic
curves $\{B_t\}$. This means that for any
$J_t$ in the path, $p'|_{B_t}$ has only
finitely many simple branch points ($z\mapsto z^2$). Then
$p'\circ\kappa$ is a DLF, which is oriented because the intersection
indices of $B$ with the fibers are positive.
\qed

\begin{remark}\rm
    1)\ \ A deeper question concerns the characterization of
    symplectic Noet\-her-Hor\-ik\-awa surfaces by invariants. It would be
    surprising, should $(c_1^2,\chi)$ alone suffice as in the
    algebraic case.\\[1.5ex]
    2)\ \ A similar argument as in the proof shows that to prove
    holomorphicity of a hyperelliptic SLF, it suffices to show that
    $(B\hookrightarrow P)$ is diffeomorphic to a holomorphic curve
    inside a rational surface. The isomorphism will automatically
    respect the fibration structures (connect the pull-back of the
    integrable complex structure to an almost complex structure $J$ that
    is compatible with $p:P\to\pr^1$ and making $B$ a $J$-holomorphic
    curve).\\[1.5ex]
    3)\ \ This proposition together with the braid monodromy computation
    for algebraic curves by specialization also shows that the isotopy
    problem for symplectic submanifolds is equivalent to the following
    group theoretic problem: The subset of half-twists in $B(S^2,d)$
    is stable under inner automorphisms. On
    $\mu$-tuples $(\tau_1,\ldots,\tau_\mu)$ of half-twists consider
    the equivalence relation generated by
    \[
      (\tau_1,\ldots,\tau_r,\tau_{r+1},\ldots,\tau_\mu)\ \sim\
      (\tau_1,\ldots,\tau_{r-1},\tau_r\tau_{r+1}\tau_r^{-1},
      \tau_r,\tau_{r+2},\ldots,\tau_\mu)
    \]
    ({\em Hurwitz-equivalence}). The claim is that
    any $(\tau_1,\ldots,\tau_\mu)$ with
    $\prod_r \tau_r\in\{1,(x_1\cdot\ldots\cdot x_{d-1})^d\}$ is
    Hurwitz-equivalent to one of
    \begin{eqnarray*}
      &&(x_1,\ldots,x_{d-1},x_{d-1},\ldots,x_1)^\frac{\mu}{2d-2}\\
      &&(x_1,\ldots,x_{d-1},x_{d-1},\ldots,x_1)^\frac{\mu-d(d-1)}{2d-2}
           (x_1,\ldots,x_{d-1})^d\\
      &&(x_1,\ldots,x_{d-2})^\frac{\mu}{d-2}\,.
    \end{eqnarray*}
    \vspace{-6ex}

\qed
\end{remark}
\pagebreak


\addcontentsline{toc}{section}{Appendix: Two-fold covers}
\noindent
{\Large\bf Appendix: Two-fold covers}
\vspace{1.5ex}

\noindent
A degree two map $\kappa:M\rightarrow N$ of $n$-dimensional, oriented,
differentiable manifolds whose set of critical values is a codimension
two, oriented submanifold $B\subset N$ is called {\em two-fold cover of
$N$ with branch locus $B$} if locally along $B$ there exist oriented
local coordinates $(x_1,\ldots,x_n)$ of $M$ and $(y_1,\ldots,y_n)$
of $N$ such that
\[
  \kappa:\ (x_1,\ldots,x_n)\ \longmapsto\
  (y_1,\ldots,y_n)\ =\ (x_1^2+x_2^2,
  x_1 x_2,x_3,\ldots,x_n)\,.
\]
The fibers having cardinality two, the corresponding unbranched covers
$M\setminus\kappa^{-1}(B)\rightarrow N\setminus B$ are Galois with group
$\gz/2\gz$. The generator of this group is the orientation preserving covering
transformation $\iota: M\rightarrow M$ swapping the two elements of a
generic fiber
(an involution).

To construct $M$ with branch locus $B\subset N$ we assume that the
homology class of $B$ is divisible by $2$ in $H_{n-2}(M;\gz)$ (we will soon
see that this
is always true for branch loci of two-fold covers).
Then $B$ is Poincar\'e-dual to the first Chern class of the square of a complex
line bundle $L$. Let $s$ be a transverse section of $L^{\otimes 2}$
with zero locus $B$. The preimage of the zero section of the quadratic map
\[
  L\ \longrightarrow\ L^{\otimes 2}\,,\quad t\longmapsto t^2-s
\]
defines an oriented, connected submanifold $M\subset L$. The bundle projection
$\kappa:L\rightarrow N$ exhibits $M$ as two-fold cover simply branched
along $B$.
We call $M=M_{L,s}\rightarrow N$ the two-fold branched cover associated
to $L$ and $s$.

It turns out that all two-fold covers arise in this way.
\begin{proposition}
  Let $\kappa: M\rightarrow N$ be a two-fold cover with branch locus $B$. Then
  there exists a complex line bundle $L$, a transverse section $s$ of
  $L^{\otimes 2}$ and a diffeomorphism of $M$ with $M_{L,s}$ over $N$.

  Moreover, the pair $(L,s)$ is uniquely determined by $M$ up to canonical
isomorphism
  (of line bundles over $N$ with section).
\end{proposition}
\pf
Away from the branch locus we define $L$ at some $x\in N$ as the space of
linear
forms $\lambda$ on ${\rm Map}(\kappa^{-1}(x);\cz)$ that change signs under
application of the involution $\iota$. Explicitely, writing
$\kappa^{-1}(x)=\{P,Q\}$
an element $f\in {\rm Map}(\kappa^{-1}(x);\cz)$ is given by the pair
$(f(P),f(Q))\in
\cz^2$ and the linear forms to be considered have the form $\lambda(a,b)=
c\cdot(a-b)$ for some $c\in\cz$. To extend $L$ across $B$ let $M\rightarrow N$
locally be given by $(z,{\bf w})\mapsto (u,{\bf v})=(z^2,{\bf w)}$
where $z,u$ are complex valued and ${\bf w}, {\bf v}$ are $(n-2)$-tuples
of coordinates. Since $\iota^*z=-z$ the linear form $\lambda_0$ sending
the values $\pm\sqrt{u}$ of the function $z$ on the fibers of
$\kappa$ to $1$ is a frame for $L$ away from $B$. The extension can
thus be defined by taking $\lambda_0$ as frame at the point of $B$ under study.

Next we observe that $L\otimes L$ has a canonical section $s$ sending
$((a,b),(c,d))$ to $(a-b)(c-d)$. This section has a simple zero along $B$
since in the local coordinates $z,{\bf w},u,{\bf v}$ we get
\[
  s(u,{\bf v})\ =\ (z-(-z))\cdot(z-(-z))\ =\ 4z^2\cdot\lambda_0^2(u,{\bf
v})\ =\
  4u\cdot\lambda_0^2(u,{\bf v})\,.
\]
To identify $\kappa:M\rightarrow N$ with $\kappa_{L,s}:
M_{L,s}\rightarrow N$ let $\pm\lambda\in L_x$ be the solutions to $t^2=s$ over
some $x\in N$. So $\kappa^{-1}(x)=\{-\lambda,\lambda\}$. Let $f(P)=1$,
$f(Q)=0$.
Then $s_x(f,f)=1$ and $\lambda(f)=1$ or $\lambda(f)=-1$. There is thus a
canonical
bijection between $\kappa^{-1}(x)$ and $\kappa_{L,s}^{-1}(x)$ by sending $P$ to
$\lambda(f)\cdot \lambda$ and $Q$ to $-\lambda(f)\cdot \lambda$. This
identification
clearly extends over $B$ as one easily checks using the local description.

In view of the functorial nature of the construction it follows also that a
diffeomorphism of two-fold coverings $M\simeq M'$ induces an isomorphism of the
associated complex line bundles $L\simeq L'$ whose square carries the canonical
sections $s,s'$ into each other.
\qed
\vspace{1ex}

\noindent
It is natural to ask to what extend the information given by $L$ and $s$ is
already
contained in $B$.
\begin{proposition}\label{twofoldcovers}
  For any 2-divisible, oriented submanifold $B\subset N$ there are exactly
  $\sharp H^1(N;\gz/2\gz)$ isomorphism classes of two-fold covers of $N$
  with branch locus $B$.

  In particular, if $H^1(N;\gz/2\gz)=0$ there
  is a one-to-one correspondence between isomorphism
  classes of two-fold covers of $N$ and 2-divisible, oriented submanifolds of
  $N$, given by the branch locus.
\end{proposition}
\pf
Let $\cala$ be the sheaf of germs of complex valued smooth functions on
$N$, and
$\cala^*$ those without zeroes (that we view as multiplicative abelian
sheaf). We will
use the following short exact sequences of sheaves
\[\begin{array}{rcccccccl}
  0&\longrightarrow&\gz/2\gz &\longrightarrow& \cala^*
  &\stackrel{\uparrow 2}{\longrightarrow}& \cala^* & \longrightarrow & 0\\[3pt]
  0& \longrightarrow & \gz &\stackrel{\cdot 2}{\longrightarrow}& \gz
  &\longrightarrow&\gz/2\gz&\longrightarrow& 0\\[3pt]
  0&\longrightarrow&\gz &\longrightarrow& \cala
  &\stackrel{\exp}{\longrightarrow}& \cala^* & \longrightarrow & 0
\end{array}\]
From the last sequence we obtain isomorphisms $H^i(N;\cala^*)\simeq
H^{i+1}(N; \gz)$ (as a soft sheaf $\cala$ has vanishing higher cohomology).
A little diagram chase shows that these are compatible with the
cohomology sequences of the other two sequences as follows:
\[\begin{array}{ccccccccc}
  H^0(N;\cala^*) &\hspace{-1ex}\stackrel{\uparrow 2}{\longrightarrow}
  \hspace{-1ex}& H^0(N;\cala^*)
  &\hspace{-1ex}\longrightarrow\hspace{-1ex}& H^1(N;\gz/2\gz)
  &\hspace{-1ex}\longrightarrow\hspace{-1ex}& H^1(N;\cala^*)
  &\stackrel{\uparrow 2}{\hspace{-1ex}\longrightarrow\hspace{-1ex}}
  &H^1(N;\cala^*)\\
  \downarrow&&  \downarrow&& \downarrow&& \downarrow&& \downarrow\\
  H^1(N;\gz) &\stackrel{\cdot 2}{\hspace{-1ex}\longrightarrow\hspace{-1ex}}&
  H^1(N;\gz) &\hspace{-1ex}\longrightarrow\hspace{-1ex}& H^1(N;\gz/2\gz)
  &\hspace{-1ex}\longrightarrow\hspace{-1ex}& H^2(N;\gz)
  &\stackrel{\cdot 2}{\hspace{-1ex}\longrightarrow\hspace{-1ex}}&H^2(N;\gz)
\end{array}\]
Here all vertical arrows are isomorphisms. Thus $H^1(N;\gz/2\gz)$ is
(non-canonically) isomorphic to the direct sum of the 2-torsion in $H^2(N;\gz)$
and the cokernel of the squaring map $\cala^*(N)\rightarrow\cala^*(N)$.
We claim that these two data act effectively on the set of isomorphism classes
of pairs $(L,s)$.

In fact, the 2-torsion of $H^2(N;\gz)$ exactly parametrizes
isomorphism classes of complex line bundles $L$ with $2 c_1(L)$ a given class
(here the Poincar\'e-dual of $B$). Once $L$ is fixed, the space of isomorphisms
$L\rightarrow L$ induces an action on the space of (transverse) sections $s$ of
$L^{\otimes 2}$. Now for any two sections $s,s'$ having the same zero locus $B$
the function $\lambda=s'/s$ has no zeros (so $\lambda\in \cala^*(N)$) and
multiplication by $\lambda$ induces the (unique) isomorphism of $L^{\otimes 2}$
carrying $s$ to $s'$. This comes from an isomorphism of $L$ iff $\lambda$
can be written
as a square. Conversely, multiplication of $s$ by any non-square $\lambda$
leads
to a pair $(L,\lambda s)$ that is not equivalent to $(L,s)$ in the
considered sense.
This concludes the proof of the proposition.
\qed
\vspace{1ex}

\begin{remark}\rm
In the holomorphic category, given $B$ and $L$ there is at most one isomorphism
class of sections $s$ of $L^{\otimes 2}$ having a holomorphic
representative. The
reason is that by the identity theorem any two holomorphic sections of a
holomorphic line bundle with the same zero divisor differ only by a constant.
Allowing changes of the complex structure of $M$ (leaving $B$ as complex
submanifold)
might however lead to different classes of sections.
\qed
\end{remark}


\vspace{0.5ex}

\noindent
{\sc \small Fakult\"at f\"ur Mathematik,
Ruhr-Universit\"at Bochum, D-44780 Bochum\\
\tt Bernd.Siebert@ruhr-uni-bochum.de\\[0.5ex]
\sc Department of Mathematics, MIT,
Cambridge, MA 02139--4307\\
\tt tian@math.mit.edu}
\end{document}